# POD-GALERKIN REDUCED ORDER METHODS FOR CFD USING FINITE VOLUME DISCRETISATION: VORTEX SHEDDING AROUND A CIRCULAR CYLINDER

GIOVANNI STABILE[1,*], SADDAM HIJAZI[1], ANDREA MOLA [1], STEFANO LORENZI[2], AND GIANLUIGI ROZZA[1]

ABSTRACT. Vortex shedding around circular cylinders is a well known and studied phenomenon that appears in many engineering fields. A Reduced Order Model (ROM) of the incompressible flow around a circular cylinder is presented in this work. The ROM is built performing a Galerkin projection of the governing equations onto a lower dimensional space. The reduced basis space is generated using a Proper Orthogonal Decomposition (POD) approach. In particular the focus is into (i) the correct reproduction of the pressure field, that in case of the vortex shedding phenomenon, is of primary importance for the calculation of the drag and lift coefficients; (ii) the projection of the Governing equations (momentum equation and Poisson equation for pressure) performed onto different reduced basis space for velocity and pressure, respectively; (iii) all the relevant modifications necessary to adapt standard finite element POD-Galerkin methods to a finite volume framework. The accuracy of the reduced order model is assessed against full order results.

1. INTRODUCTION

A large part of physical systems is described by partial differential equations and their numerical solution is essential in many engineering fields. Even though several progresses have been made over the last decades, the numerical solution of fluid dynamics problems, using standard finite element methods (FEM), spectral element methods (SEM), finite volume methods (FVM) or finite differences methods (FDM), may be extremely expensive by a computational standpoint. The development of efficient and reliable Reduced Order Models (ROMs) could be a great advantage especially when dealing with control, optimization and uncertainty quantification problems, where a large number of different system configurations are in need of being tested.

Vortex Induced Vibrations (VIVs) are an important phenomena in many different engineering fields where it has been observed and studied either both for air or water flows for many years [27, 28]. The importance of studying such a problem comes from the fact that it can be the source of evident damage or failure of the engineering system. Disregarding VIVs through the design process can in fact lead to severe structural failures. This phenomenon is caused by the oscillating flow arising from the alternate vortex shedding. Among all possible existing phenomena that may occur on flexible cylindrical structures, VIVs are potentially one of the most dangerous and hard to predict. If a rigid and fixed cylinder is considered, the frequency of the vortex shedding phenomenon $f_v$ can be deduced with $f_v = StU/D$ [40] in which $St$ is the *Strouhal* number, $U$ is the free stream velocity and $D$ is the diameter of the cylinder. For the particular case of the vortex shedding phenomenon around a circular cylinder [45, 44] one can find several attempts to create a ROM. In some cases no attempt to model completely the flow field is done such as in [15, 13, 38] where only the lift and drag forces acting on the cylinder are modelled using wake oscillator models. When the interest is into the complete reconstruction of the flow field, Reduced Basis (RB) method can be applied. In this method the governing equations, describing the phenomenon, are projected onto a low dimensional space called the reduced basis space [16, 30] that is optimally constructed starting from high fidelity simulations. In particular in this work the ROM is constructed using a POD-Galerkin approach [25, 1, 21, 5, 20, 6, 2]. As previously mentioned, since the main purpose of this paper is the correct modelling of the vortex shedding phenomenon, particular attention is paid to the reconstruction in the reduced order model of the pressure field. Differently from [21], in this paper we propose a ROM that includes also the Poisson equation for the pressure modelling. To the best of the authors knowledge, this is the first attempt to use this approach, commonly employed in the FE framework, in the POD-Galerkin approach using finite volume discretization.

[1]SISSA, International School for Advanced Studies, Mathematics Area, mathLab Trieste, Italy.
[2]Department of Energy, Politecnico di Milano, Italy.
*E-mail addresses*: gstabile@sissa.it, shijazi@sissa.it, amola@sissa.it, stefano.lorenzi@polimi.it , grozza@sissa.it.
2010 *Mathematics Subject Classification.* 78M34, 97N40, 35Q35.
[*]Corresponding Author.

The work is organized as follows. In § 2 the governing equations of the physical model and the high fidelity (HF) discretisation techniques used to solve the full order model are presented. The development of the ROM is introduced in § 3 and in § 4 the numerical example, regarding the vortex shedding phenomenon around a circular cylinder, is analysed. Finally in § 5 conclusions and suggestions for future developments are given.

## 2. The Full Order Model

The physical model is described below by using the parametrized incompressible unsteady Navier-Stokes equations. They consists into the well known conservation of momentum law and continuity equations. In an Eulerian framework they are expressed by:

$$
(2.1) \quad \begin{cases} \frac{\partial \boldsymbol{u}}{\partial t} + (\boldsymbol{u} \cdot \boldsymbol{\nabla})\boldsymbol{u} - \boldsymbol{\nabla} \cdot \nu \boldsymbol{\nabla} \boldsymbol{u} = -\boldsymbol{\nabla} p & \text{in } \Omega_f \times [0,T], \\ \boldsymbol{\nabla} \cdot \boldsymbol{u} = \boldsymbol{0} & \text{in } \Omega_f \times [0,T], \\ \boldsymbol{u}(t,x) = \boldsymbol{f}(x,\mu) & \text{on } \Gamma_{In} \times [0,T], \\ \boldsymbol{u}(t,x) = \boldsymbol{0} & \text{on } \Gamma_0 \times [0,T], \\ (\nu \boldsymbol{\nabla} \boldsymbol{u} - p\boldsymbol{I})\boldsymbol{n} = \boldsymbol{0} & \text{on } \Gamma_{Out} \times [0,T], \\ \boldsymbol{u}(0,\boldsymbol{x}) = \boldsymbol{k}(\boldsymbol{x}) & \text{in } (\Omega_f, 0) \end{cases}
$$

where $\Gamma = \Gamma_{In} \cup \Gamma_0 \cup \Gamma_{Out}$ is the boundary of the fluid domain $\Omega_f$ and is composed by three different parts $\Gamma_{In}$, $\Gamma_{Out}$ and $\Gamma_0$ that indicate respectively inlet boundary, outlet boundary and physical walls. $\boldsymbol{u}$ is the flow velocity vector, $t$ is the time, $\nu$ is the fluid kinematic viscosity, and $p$ is the normalized pressure, which is divided by the fluid density $\rho_f$, $\boldsymbol{f}$ is a generic function that gives the value of the velocity on the inlet $\Gamma_{In}$ and it is parametrised through the scalar quantity $\mu$. $\boldsymbol{k}$ is the initial velocity field and $T$ is the time window we considered. Since in the present work the problem is solved using a finite volume discretisation technique [41, 24], where the standard is to work with a Poisson equation for pressure rather than directly with the continuity equation, the system of equations 2.1 is modified into:

(2.2)
$$\begin{cases} \frac{\partial \boldsymbol{u}}{\partial t} + (\boldsymbol{u} \cdot \boldsymbol{\nabla})\boldsymbol{u} - \boldsymbol{\nabla} \cdot \nu \boldsymbol{\nabla} \boldsymbol{u} = -\boldsymbol{\nabla} p & \text{in } \Omega_f \times [0, T], \\ \boldsymbol{u}(t, x) = \boldsymbol{f}(x, \mu) & \text{on } \Gamma_{In} \times [0, T], \\ \boldsymbol{u}(t, x) = \boldsymbol{0} & \text{on } \Gamma_0 \times [0, T], \\ \nabla \boldsymbol{u} \cdot \boldsymbol{n} = \boldsymbol{0} & \text{on } \Gamma_{Out} \times [0, T], \\ \boldsymbol{u}(0, \boldsymbol{x}) = \boldsymbol{k}(\boldsymbol{x}) & \text{in } (\Omega_f, 0), \\ + \\ \Delta p = -\nabla \cdot (\boldsymbol{u} \cdot \nabla)\boldsymbol{u} & \text{in } \Omega_f \times [0, T], \\ \nabla p \cdot \boldsymbol{n} = 0 & \text{on } \Gamma \setminus \Gamma_{Out} \times [0, T], \\ p = 0 & \text{on } \Gamma_{Out} \times [0, T]. \end{cases}$$

In the above system of equations all the quantities assume the same meaning of those presented in System 2.1 and the Poisson equation for pressure is obtained taking the divergence of the momentum equation and exploiting the divergence free costraint. The two equations are solved using a segregated approach and more details are given in § 2.1. Historically the FVM discretisation technique is widely used in industrial applications and for flows characterized by high values of the *Reynolds* numbers. An advantage of the FVM is that of ensuring local enforcement of the conservative law since equations are written in conservative form. In this work the high fidelity simulations are carried out making use of the finite volume C++ library OpenFOAM® (OF) [42].

2.1. **The Finite Volume discretisation.**
As mentioned in § 1 the governing equations are discretised in space using a finite volume approximation. Once a suitable polygonal tessellation is chosen, the system of partial differential equations in (2.1) is written in integral form over a control volume. In the present two-dimensional framework the tessellation is represented by a subdivision of the domain into a finite number of non-overlapping polygonal cells. The dimension of the full order model, which consists into the number of degrees of freedom of the discretised problem, will be henceforth indicated with $N_h$. The strategies for the discretisation of both the momentum and continuity equation are briefly reported in the following.

2.1.1. *Momentum equation.*
The momentum balance equation is written for each volume $V_i$ in integral form as:

(2.3) $$\int_{V_i} \frac{\partial}{\partial t}\boldsymbol{u_t} dV + \int_{V_i} (\boldsymbol{u} \cdot \nabla)\boldsymbol{u} dV - \int_{V_i} \nabla \cdot \nu \nabla \boldsymbol{u} dV + \int_{V_i} \nabla p dV = 0.$$

The gradient terms and in particular the gradient of pressure, making use of the Gauss's theorem, are discretised as:

(2.4) $$\int_{V_i} \nabla p dV = \int_{S_i} d\boldsymbol{S} \cdot p \approx \sum_f \boldsymbol{S_f} p_f,$$

where $\boldsymbol{S_f}$ is the area vector of each face of the control volume and $p_f$ is the value of pressure at the center of the faces (Figure 1).

Making use of the Gauss's theorem, the convective term is discretised as follow:

(2.5) $$\int_{V_i} (\boldsymbol{u} \cdot \nabla)\boldsymbol{u} dV = \int_{S_i} (d\boldsymbol{S} \cdot \boldsymbol{u_f})\boldsymbol{u_f} \approx \sum_f (\boldsymbol{S_f} \cdot \boldsymbol{u_f})\boldsymbol{u_f} = \sum_f F_f \boldsymbol{u_f},$$

$\boldsymbol{u_f}$ is the velocity vector evaluated at the center of each face of the control volume. Since the unknowns of the mathematical problem are the velocity values at the cell centers, $u_f$ must be be obtained from the cell center values with suitable interpolation schemes. Possible alternatives could be a central, upwind, second order upwind and blended differencing schemes. In this particular case a second order upwind differencing scheme is used for the convective term. $F_f = \bm{S_f} \cdot \bm{u}_f$ is the mass flux through each face of the control volume. The diffusion term is discretised as:

$$\int_{V_i} \nabla \cdot \nu \nabla \bm{u} \, dV = \int_{S_i} d\bm{S} \cdot \nu \nabla \bm{u} \approx \sum_f \nu \bm{S_f} \cdot (\nabla \bm{u})_f, \tag{2.6}$$

where $(\nabla \bm{u})_f$ is the gradient of $\bm{u}$ at the faces. In case the value of the gradient $\nabla \phi_i$ of a generic conservative variable $\phi$ at the center of the cell is needed, as in equation 2.6, this can be computed dividing the expression in equation 2.4 by the volume of the cell $V_i$. In case of orthogonal meshes (i.e. the face dividing two cells is orthogonal with respect to the distance connecting the two cell centers) the term $\bm{S_f} \cdot (\nabla \bm{u})_f$ could be computed using:

$$\bm{S_f} \cdot (\nabla \bm{u})_f = |\bm{S_f}| \frac{\bm{u}_N - \bm{u}_P}{|\bm{d}|}, \tag{2.7}$$

where $\bm{u}_N$ and $\bm{u}_P$ are the velocities at the centers of two neighboring cells and $\bm{d}$ is the distance vector connecting the two cell centers (see Figure 1). In the case of non-orthogonal meshes one needs to correct the above scheme. In this work this term is split into two contributions, an orthogonal one and a non-orthogonal one [18]:

$$\bm{S_f} \cdot (\nabla \bm{u})_f = |\bm{\Delta}| \frac{\bm{u}_N - \bm{u}_P}{|\bm{d}|} + \bm{k} \cdot (\nabla \bm{u})_f, \tag{2.8}$$

where the two vectors $\bm{\Delta}$ and $\bm{k}$ satisfy $\bm{S_f} = \bm{\Delta} + \bm{k}$. The first vector $\bm{\Delta}$ is chosen parallel to $\bm{S_f}$. The term $(\nabla \bm{u})_f$ is obtained through interpolation of the the values of the gradient at the cell centers $(\nabla \bm{u})_N$ and $(\nabla \bm{u})_P$ in which the subscripts $N$ and $P$ (Figure 1) indicate the values at center of the cells of the two neighboring cells. There are different strategies to determine the vectors $\bm{\Delta}$ and $\bm{k}$ such as minimum correction approach, orthogonal correction approach and over-relaxed approach, more details can be found in [18].

2.1.2. *Poisson equation for pressure.*
As mentioned in § 1 the coupling between momentum conservation and the Poisson equation for pressure is treated making use of a segregated approach. The PIMPLE algorithm [24] is used, it consists into a combination of an inner correction cycle using the PISO [17] algorithm and an outer correction procedure performed using the SIMPLE [29] algorithm. The two considered equations are the momentum equation and the Poisson equation for pressure, which as already mentioned, is obtained taking the divergence of the momentum equation and exploiting the continuity equation:

$$\Delta p = -\nabla \cdot (\bm{u} \cdot \nabla) \bm{u}. \tag{2.9}$$

The use of a Poisson equation for pressure is considered also at reduced order level where the two equations are instead solved in a monolithic way.

## 3. The Reduced Order Model

This section introduces the derivation of the ROM and the necessary modifications in order to be able to adapt a standard FEM-Galerkin ROM to a FVM-Galerkin framework. Here few details are only recalled, for further details one may see [21]. The relevant details introduced here are:

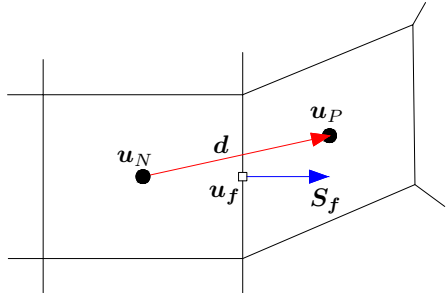

FIGURE 1. Sketch of a finite volume in 2 dimensions

- Introduction of a reduced basis space also for the mass flux term of equation (2.5) with a similar approach presented also in [21]. This reduced basis space is the one used during the projection phase of the momentum equation.
- Differently to what proposed in [21, 5], where only the momentum equation is considered and where it is assumed that velocity and pressure share the same temporal coefficients, the Poisson equation for pressure reported in Equation (2.9) is projected onto the POD pressure modes in order to enforce the continuity equation constraint. Such approach, that is proposed in literature by several authors [1, 7] for FEM approximations, is here adapted to a FVM framework.

The reduced order model is obtained performing a Galerkin projection onto the space spanned by the reduced basis modes and approximating the fields with the following expansions:

$$(3.1) \qquad \begin{pmatrix} \boldsymbol{u}(\boldsymbol{x},t) \\ F(\boldsymbol{x},t) \end{pmatrix} \approx \begin{pmatrix} \boldsymbol{u}_r(\boldsymbol{x},t) \\ F_r(\boldsymbol{x},t) \end{pmatrix} = \sum_{i=1}^{N_u} a_i(t) \begin{pmatrix} \boldsymbol{\varphi}_i(\boldsymbol{x}) \\ \psi_i(\boldsymbol{x}) \end{pmatrix},$$

$$(3.2) \qquad p(\boldsymbol{x},t) \approx p_r(\boldsymbol{x},t) = \sum_{i=1}^{N_p} b_i(t) \chi_i(\boldsymbol{x}),$$

where the $a_i$ and $b_i$ are temporal coefficients and $\boldsymbol{\varphi_i}$, $\psi_i$ and $\chi_i$ are the modes of the reduced basis spaces for velocity, mass flux and pressure, respectively. $N_u$ and $N_p$ define the dimension of the reduced basis spaces for velocity/mass flux and pressure. Clearly, $N_u$ and $N_p$ are not constrained to have coincident values. In the above expansions the following assumptions are considered:

- Velocity $\boldsymbol{u}$ and mass flux $F$ fields are approximated using the same temporal coefficients. This assumption is reasonable since also the mass fluxes over the surfaces of the finite volumes depend strongly on the velocity itself. The mass flux term which is a scalar field defined on the surfaces of all the cells, as indicated in equation 2.5, is defined as a product of the velocity at face, which is obtained through interpolation, and the area of the face. Moreover the reduced basis space of mass flux and velocity respectively have the same dimension $N_u$.
- Pressure field $p$ is approximated using different temporal coefficients respect to the velocity/mass flux fields and for this reason during the projection phase both the momentum conservation and Poisson equation for pressure must be considered. This space can have also a different dimension respect to the reduced basis space considered for velocity.

3.1. **Generation of the POD spaces.**
In order to create a reduced basis space onto which the governing equations are projected, one can find many techniques in literature such as the Proper Orthogonal Decomposition (POD), The Proper Generalized Decomposition (PGD), as well as Reduced Basis (RB) method with a greedy approach. For more details about the different methods the reader may see [33, 9, 16, 30, 10, 12]. In this work the POD approach is used. The POD consists into the decomposition of the flow fields into temporal coefficients $a_i(t)$ and orthonormal spatial bases $\boldsymbol{\varphi}_i(\boldsymbol{x})$:

$$(3.3) \qquad \boldsymbol{u}(\boldsymbol{x}, t) = \sum_{i=1}^{N_s} a_i(t) \boldsymbol{\varphi}_i(\boldsymbol{x}),$$

where $\boldsymbol{\varphi}_i(\boldsymbol{x})$ are orthonormal spatial bases that minimizes the average of the error between the snapshots and their orthogonal projection onto the bases and $N_s$ is the number of considered snapshots. The POD space $\mathbb{V}_{POD} = \text{span}(\boldsymbol{\varphi_1}, \boldsymbol{\varphi_2}, \ldots, \boldsymbol{\varphi_{N_s}})$ is then constructed solving the following minimization problem:

$$(3.4) \qquad \mathbb{V}_{POD} = \arg\min \frac{1}{N_s} \sum_{n=1}^{N_s} \| \boldsymbol{u}_n - \sum_{n=1}^{N_s} (\boldsymbol{u}_n, \boldsymbol{\varphi}_i)_{L^2(\Omega)} \boldsymbol{\varphi}_i \|_{L^2(\Omega)}^2,$$

where $\boldsymbol{u}_n$ is a general snapshot of the velocity field at time $t = t_n$. This problem can be solved computing a singular value decomposition $\boldsymbol{\mathcal{U}} = \boldsymbol{\mathcal{W}^u \Sigma^u \mathcal{V}^{uT}}$ of the so called snapshots matrix $\boldsymbol{\mathcal{U}} \in \mathbb{R}^{N_h \times N_s}$. Where $\boldsymbol{\mathcal{U}} = [\boldsymbol{u_1}, \boldsymbol{u_2}, \ldots, \boldsymbol{u_{N_s}}]$ contains the flow fields for all the different time steps, $\boldsymbol{\mathcal{W}^u} \in \mathbb{R}^{N_h \times N_h}$ is a rectangular matrix of left singular vectors, $\boldsymbol{\mathcal{V}^u} \in \mathbb{R}^{N_s \times N_s}$ is a square matrix of right singular vectors, and $\boldsymbol{\Sigma^u} \in \mathbb{R}^{N_h \times N_s}$ is a diagonal matrix of eigenvalues. The POD modes $\boldsymbol{\varphi_i}$ are then given by the columns of the matrix $\boldsymbol{\mathcal{W}^u}$. This approach might be however computationally expensive, especially increasing the dimension of the grid used to discretise the domain. An equivalent and more efficient way to tackle this problem, based on the method of snapshots, firstly introduced in [37], consists in solving the eigenvalue problem:

$$(3.5) \qquad \boldsymbol{C^u Q^u} = \boldsymbol{Q^u \lambda}^u,$$

where $\boldsymbol{C^u} \in \mathbb{R}^{N_s \times N_s}$ is the correlation matrix of the velocity field snapshots, $\boldsymbol{Q^u} \in \mathbb{R}^{N_s \times N_s}$ is a square matrix whose columns are the eigenvectors and $\boldsymbol{\lambda^u} \in \mathbb{R}^{N_s \times N_s}$ is a diagonal matrix containing the eigenvalues $\lambda_{ii}^u$. The correlation matrix can be determined using:

$$(3.6) \qquad C_{ij}^u = (\boldsymbol{u_i}, \boldsymbol{u_j})_{L^2(\Omega)}$$

where $(\cdot, \cdot)_{L^2(\Omega)}$ is the $L^2$ inner product over the domain $\Omega$.

*Remark* 3.1. Normally, in a standard finite element framework, the $H^1$ norm is preferred for the velocity field since its natural functional space is $H^1(\Omega)$. Here it is decided to use the $L^2$ norm for both the pressure and the velocity fields. In a finite volume setting in fact both velocity and pressure belong to discontinuous spaces and, as illustrated in equation 2.4, in order to compute the gradients necessary for the $H^1$ norm evaluation one would introduce further discretization error.

The POD modes can be finally obtained with:

$$(3.7) \qquad \boldsymbol{\varphi_i} = \frac{1}{\sqrt{\lambda_{ii}^u}} \boldsymbol{\mathcal{U} Q_i^u}.$$

The same procedure can be repeated also for the pressure field considering the snapshots matrix $\boldsymbol{\mathcal{P}} = [p_1, p_2, \ldots, p_{N_s}]$ one can compute the correlation matrix of the pressure field snapshots $C^p$ and solve a similar eigenvalue problem $C^p Q^p = Q^p \lambda^p$. The POD modes $\chi_i$ for the pressure field can be computed with:

$$\chi_i = \frac{1}{\sqrt{\lambda_{ii}^p}} \mathcal{P} Q_i^p. \tag{3.8}$$

The modes for mass flux field are obtained using the same eigenvectors and eigenvalues computed solving the eigenvalue problem of the velocity field and are expressed by:

$$\psi_i = \frac{1}{\sqrt{\lambda_{ii}^u}} \mathcal{F} Q_i^u, \tag{3.9}$$

where, again, $\mathcal{F} = [F_1, F_2, \ldots, F_{N_s}]$ is a snapshots matrix containing the mass flux field at different time steps. For what concern the basis for the mass flux term it is decided to use the same eigenvector of the eigenvalue problem solved for the velocity field because it is assumed that mass flux and velocity share the same temporal coefficients. More details are given in the next subsection.

### 3.2. Galerkin projection onto the POD space.

In this section the Galerkin projection of the governing equations onto the POD space is highlighted and discussed. The idea here is to exploit both the momentum conservation and continuity equation.

#### 3.2.1. ROM for velocity - Momentum equation.

The reduced order model of the momentum equation is obtained performing an $L^2$ orthogonal projection onto the reduced bases space $\mathbb{V}_{POD}$ spanned by the POD velocity modes with a procedure similar to what presented in [1].

$$(\boldsymbol{\varphi_i}, \boldsymbol{u_t} + (\boldsymbol{u} \cdot \nabla)\boldsymbol{u} - \nu \Delta \boldsymbol{u} + \nabla p)_{L^2(\Omega)} = 0. \tag{3.10}$$

Respect to what presented in [1] here also the gradient of pressure is considered inside the momentum equation. This term is considered also in [5, 21] but, there, it is assumed that velocity and pressure modes share the same temporal coefficients. More details about the treatment of this term are given in § 3.2.2. Substituting the POD approximations of $\boldsymbol{u}$, $F$ and $p$ into equation (3.10) and exploiting the orthogonality of the POD modes $\boldsymbol{\varphi_i}$ one obtains the following dynamical system:

$$\dot{\boldsymbol{a}} = \nu \boldsymbol{B} \boldsymbol{a} - \boldsymbol{a}^T \boldsymbol{C} \boldsymbol{a} - \boldsymbol{K} \boldsymbol{b}, \tag{3.11}$$

where $\boldsymbol{a}$ and $\boldsymbol{b}$ are vectors containing all the temporal coefficients $a_i(t)$ and $b_i(t)$ and the terms inside equation (3.11) read:

$$B_{ij} = (\boldsymbol{\varphi_i}, \Delta \boldsymbol{\varphi_j})_{L^2(\Omega)}, \tag{3.12}$$

$$C_{ijk} = (\boldsymbol{\varphi_i}, \nabla \cdot (\psi_j, \boldsymbol{\varphi_k}))_{L^2(\Omega)}, \tag{3.13}$$

$$K_{ij} = (\boldsymbol{\varphi_i}, \nabla p_j)_{L^2(\Omega)}. \tag{3.14}$$

In Equation (3.13) the term $\nabla \cdot (\psi_j, \boldsymbol{\varphi_k})$, with an abuse of notation, is used to indicate the convective term. This term is obtained exploiting equation (2.5) and the velocity approximation:

$$\int_{V_i} (\boldsymbol{\varphi_j} \cdot \nabla) \boldsymbol{\varphi_k} \mathrm{d}V = \int_{\partial V_i} (\boldsymbol{dS} \cdot \boldsymbol{\varphi_j}) \boldsymbol{\varphi_k} = \sum_f \boldsymbol{S_f} \cdot \boldsymbol{\varphi_{j,f}} \boldsymbol{\varphi_{k,f}} = \sum_f \psi_{j,f} \boldsymbol{\varphi_{k,f}} \tag{3.15}$$

For each mode of the velocity and mass fluxes POD spaces, the velocity base $\boldsymbol{\varphi_k}$, which is defined at the center of each cell, must be interpolated in order to obtain the value $\boldsymbol{\varphi_{k,f}}$ at the center of each face. Multypling the value of the velocity base at the faces $\boldsymbol{\varphi_{k,f}}$ by the base for mass flux $\psi_j$, which is already defined at the center of each face, and performing a summation over the faces of each face, it is possible to evaluate the convective term at the center of each cell.

*Remark* 3.2. To ensure an efficient offline–online decomposition, even though we are dealing with an affine parametric dependence problem, in the case of the non-linear convective term, further difficulties arise. The third order tensor $C_{ijk}$ is stored [31, 34] to deal with the non-linear term in the present work. During the online solution, at each fixed point iteration of the solution procedure the $i^{\text{th}}$ component of the the residual due to the non-linear term is evaluated as:

$$(3.16) \qquad R_i = \boldsymbol{a}^T \boldsymbol{C_{i\bullet\bullet}} \boldsymbol{a}.$$

Since the dimension of the $\boldsymbol{C}$ tensor is growing with the cube of the number of basis functions employed for the velocity space, this approach may lead in some cases to high storage costs. In all the test cases presented in this work, the small dimension ($N < 10$) of the reduced space did not lead to such problem. Yet, if richer reduced spaces are used, possible alternatives to the present approach could be using EIM-DEIM [43, 4] approaches or Gappy-POD [8].

### 3.2.2. *ROM for pressure - Poisson equation for pressure.*

System (3.11) accounts for $N_u$ equations, given by the momentum equation projected on each of the velocity modes. Yet, the system presents $N_u + N_p$ unknowns given by the temporal coefficients of velocity $\boldsymbol{a}$ and pressure $\boldsymbol{b}$. Additional equations are required to close the problem. In the reduced framework, the continuity equation cannot be directly exploited because the velocity modes, which are generated with divergence free snapshots, are in turn divergence free up to numerical precision. The additional unknowns inside (3.11) are multiplied by the gradient of pressure that in many cases is neglected [26, 7]; in fact, in many contributions available in literature no attempt to recover the pressure term is performed. The projection of the pressure gradient onto the POD spaces is in fact zero for the case of enclosed flows as presented in [11, 22, 25] or in the case of inlet-outlet problems with outlet far from the obstacle [1]. However in many applications the pressure term is needed as highlighted in [26] and cannot be neglected. In the analysed case moreover, since the interest is into the reconstruction of the fluid forces acting onto the cylinder surface, the reconstruction of the reduced pressure term is crucial.

According to [7] in literature one can find basically two different approaches for pressure ROMs depending if they use pressure POD modes or not. In methods using only a POD basis for velocity the momentum equation without the gradient of pressure term is solved. The pressure field is then a posteriori reconstructed exploiting the Poisson equation for pressure 2.9. The velocity field on the right hand side of the equation is approximated with the reduced order model approximation of the velocity:

$$(3.17) \qquad -\Delta p_r = \nabla \cdot \left( \left( \sum_{i=1}^{N_u} a_i \boldsymbol{\varphi}_i \cdot \nabla \right) \sum_{j=1}^{N_u} a_j \boldsymbol{\varphi}_j \right) = \sum_{i=1}^{N_u} \sum_{j=1}^{N_u} a_i a_j \nabla \cdot \left( \left( \boldsymbol{\varphi_i} \cdot \nabla \right) \boldsymbol{\varphi_j} \right),$$

Since the temporal coefficients do not depend on space, the pressure term can be recovered with:

$$(3.18) \qquad p_r = \sum_{i=1}^{N_u} \sum_{j=1}^{N_u} a_i a_j p_{0ij},$$

$$(3.19) \qquad -\Delta p_{0ij} = \nabla \cdot \left( \left( \boldsymbol{\varphi_i} \cdot \nabla \right) \boldsymbol{\varphi_j} \right)$$

Using such an approach the term $p_{0ij}$ of equation 3.19 needs to be precomputed and this implies the resolution of $(N_u + 1)N_u/2$ Poisson problems (due to the symmetry of $p_{0ij}$) that have the dimension of the full order problem.

In ROMs exploting also the pressure modes two different approaches can be found whether they use only the momentum equation or also the continuity/Poisson equation.

In the first approach [21, 5] it is assumed that velocity and pressure share the same temporal coefficients ($\boldsymbol{a} = \boldsymbol{b}$) and only the momentum equation is exploited at reduced order level.

In the second approach it is assumed that velocity and pressure at reduced order level are approximated with different temporal coefficient ($\boldsymbol{a} \neq \boldsymbol{b}$) and also the continuity equation or the Poisson equation for pressure are exploited. Among the methods using also pressure modes, some work directly on the system composed by the continuity and momentum equation and use ad hoc stabilization techniques to enforce the well posedness of the problem [3, 35, 7, 32]; other methods work with the momentum equation, without the gradient of pressure term, and the Poisson equation for pressure [1, 7]. In the latter approaches, the momentum equation is decoupled from the Poisson equation for pressure and the pressure is reconstructed in a post processing stage after the resolution of the momentum Poisson equation.

Here the second approach is used, two different coefficients depending on time are considered (one for velocity and the other one for pressure) and then the Poisson equation for pressure is exploited. Respect to what done in [1, 7] as shown in equation 3.10, the gradient of pressure term is not neglected in the momentum equation and this gives rise to a coupled system also at reduced order level.

We remark that in the finite volume solver employed, the PIMPLE algorithm for pressure coupling is indeed based on the Poisson equation for pressure. Thus, the ROM procedure used has also the benefit of making the ROM equations consistent with the high fidelity model ones.

Poisson equation is projected onto the POD space spanned by the pressure modes $\chi_i$ and after integration by part of the Laplacian term inside equation (2.9) one obtains:

$$(3.20) \qquad (\nabla \chi_i, \nabla p)_{L^2(\Omega)} = (\chi_i, \nabla \cdot ((\boldsymbol{u} \cdot \nabla)\boldsymbol{u}))_{L^2(\Omega)},$$

which can be rewritten in matrix form as:

$$(3.21) \qquad \boldsymbol{D}\boldsymbol{b} = -\boldsymbol{E} + \boldsymbol{a}^T \boldsymbol{G} \boldsymbol{a},$$

and the terms inside (3.21) read:

$$(3.22) \qquad D_{ij} = (\nabla \chi_i, \nabla \chi_j)_{L^2(\Omega)},$$
$$(3.23) \qquad E_i = (\nabla \chi_i, \nabla \bar{p})_{L^2(\Omega)},$$
$$(3.24) \qquad G_{ijk} = (\chi_i, \nabla \cdot (\nabla \cdot (\psi_j, \boldsymbol{\varphi_k})))_{L^2(\Omega)}.$$

In the above expressions we made use of the assumption that the pressure term is decomposed into a mean and a fluctuating term $p = \bar{p} + p'$. The bases for pressure $\chi_i$ are then constructed starting from the snapshots matrix of the fluctuating pressure $\boldsymbol{\mathcal{P}'} = [p'_1, p'_2, \ldots, p'_{N_s}]$. The pressure field is then approximated with:

$$(3.25) \qquad p(\boldsymbol{x}, t) \approx p_r(\boldsymbol{x}, t) = \bar{p} + \sum_{i=1}^{N_p} b_i(t) \chi_i(\boldsymbol{x}).$$

The approach for the treatment of the non-linear term $\boldsymbol{G}$ is analogous to what done for the convection term $\boldsymbol{C}$ of the momentum equation. Also the mass flux bases are considered as a consequence of the finite volume discretisation and a third order tensor is stored. During the online procedure it is exploited the same approach used in equation 3.16.

3.3. **The boundary conditions.**
The interest here is to deal with parametrized boundary conditions also at reduced order level. In literature different approaches to enforce the BCs in the ROM can be found. In this section it is explained how the Dirichlet boundary conditions are enforced at the reduced order level. The penalty method is used in [21, 19, 36] where the BCs are imposed weakly using a penalty term, however, this method relies on a penalty parameter that has to be tuned with a sensitivity analysis [36]. In this work, in order to enforce the BCs at the ROM level a control function method is used. Within this method, before applying the POD, the inhomogenous boundary conditions are removed from the original snapshots. Using such an approach it is possible to produce homogeneous basis functions and later on, at reduced order level, is possible to deal with any boundary condition (of course it must be sufficiently close to those used to train the ROM model). The problem is then solved and the lifting function is added again to the solution. Only boundary conditions that can be parametrized with a single time-dependent coefficient as in Graham [14] are considered. To retain the divergence-free property of the snapshots the lifting function has also to be divergence free. In particular it is chosen to use as lifting function the arithmetic average of the velocity snapshots, that is opportunely scaled in order to have the desired value at the Dirichlet boundary. Each snapshot of the velocity snapshots matrix is then modified as:

$$\boldsymbol{u}'_i = \boldsymbol{u}_i - u_D(t)\boldsymbol{\phi_c}, \tag{3.26}$$

where $\boldsymbol{\phi_c}$ is a function that has unitary value at the reference point chosen for the scaling Dirichlet boundary. This lifting function can be evaluated as the arithmetic average of the velocity snapshots $\boldsymbol{u_m}$ opportunely divided by its own value $u_{m,r}$ at the reference point on the Dirichlet boundary:

$$\boldsymbol{\phi_c} = \frac{1}{N_s u_{m,r}} \sum_{i=1}^{N_s} \boldsymbol{u}_i. \tag{3.27}$$

The POD is then applied to the snapshots matrix $\boldsymbol{\mathcal{U}}' = [\boldsymbol{u}'_1, \boldsymbol{u}'_2, \ldots, \boldsymbol{u}'_{N_s}]$ that contains only snapshots with homogeneous boundary conditions. The velocity field is then approximated as:

$$\boldsymbol{u}(\boldsymbol{x}, t) \approx u_D(t)\boldsymbol{\phi_c} + \sum_{i=1}^{N_u} a_i(t)\boldsymbol{\varphi}_i(\boldsymbol{x}), \tag{3.28}$$

where $u_D(t)$ is a scaling factor depending on time that assumes the value of the Dirichlet BC at the reference point. For sake of simplicity, since time dependent BCs are not considered, the time dependency on $u_D(t)$ will be henceforth omitted. The same procedure is also repeated for the mass fluxes where the term $F_c$ is the mass flux associated with the velocity field $\boldsymbol{\phi_c}$ and the mass flux is then approximated as:

$$F(\boldsymbol{x}, t) \approx u_D F_c + \sum_{i=1}^{N_u} a_i(t)\psi_i(\boldsymbol{x}). \tag{3.29}$$

During the projection stage illustrated in § 3.2.1 and § 3.2.2 also the above modified approximations of the velocity and mass flux fields have to be considered. The Galerkin

VORTEX SHEDDING AROUND A CIRCULAR CYLINDER USING A POD-GALERKIN METHOD11...true

projection produces then some additional terms inside the coupled dynamical system that now reads:

$$\begin{cases} \dot{\boldsymbol{a}} = \boldsymbol{A}_{BC} + (\boldsymbol{B} + \boldsymbol{B}_{BC})\boldsymbol{a} - \boldsymbol{a}^T\boldsymbol{C}\boldsymbol{a} - \boldsymbol{K}\boldsymbol{b} \\ \boldsymbol{b} = \boldsymbol{D}^{-1}(\boldsymbol{E} + \boldsymbol{E}_{BC} + \boldsymbol{F}_{BC}\boldsymbol{a} + \boldsymbol{a}^T\boldsymbol{G_1}\boldsymbol{a}), \end{cases} \quad (3.30)$$

where $\boldsymbol{A}_{BC}$, $\boldsymbol{B}_{BC}$, $\boldsymbol{E}_{BC}$ and $\boldsymbol{F}_{BC}$ are equal to

$$\boldsymbol{A}_{BC} = \nu u_D \boldsymbol{A_1} - u_D^2 \boldsymbol{A_2}, \quad (3.31)$$

$$\boldsymbol{B}_{BC} = -u_D \boldsymbol{B_1} - u_D \boldsymbol{B_2}, \quad (3.32)$$

$$\boldsymbol{E}_{BC} = u_D^2 \boldsymbol{E_1}, \quad (3.33)$$

$$\boldsymbol{F}_{BC} = u_D(\boldsymbol{F_1} + \boldsymbol{F_2}). \quad (3.34)$$

The terms are obtained through Galerkin projection of the momentum equation and the Poisson equation for pressure onto the POD velocity and pressure space, respectively:

$$A_{1i} = (\boldsymbol{\varphi_i}, \Delta \boldsymbol{\phi_c})_{L^2(\Omega)}, \quad (3.35)$$

$$A_{2i} = (\boldsymbol{\varphi_i}, \nabla \cdot (F_c, \boldsymbol{\phi_c}))_{L^2(\Omega)}, \quad (3.36)$$

$$B_{1ij} = (\boldsymbol{\varphi_i}, \nabla \cdot (\psi_j, \boldsymbol{\phi_c}))_{L^2(\Omega)}, \quad (3.37)$$

$$B_{2ij} = (\boldsymbol{\varphi_i}, \nabla \cdot (F_c, \boldsymbol{\varphi_j}))_{L^2(\Omega)}, \quad (3.38)$$

$$E_{1i} = (\chi_i, \nabla \cdot (\nabla \cdot (F_c, \boldsymbol{\phi_c})))_{L^2(\Omega)}, \quad (3.39)$$

$$F_{1ij} = (\chi_i, \nabla \cdot (\nabla \cdot (\psi_j, \boldsymbol{\phi_c})))_{L^2(\Omega)}, \quad (3.40)$$

$$F_{2ij} = (\chi_i, \nabla \cdot (\nabla \cdot (F_c, \boldsymbol{\varphi_j})))_{L^2(\Omega)}. \quad (3.41)$$

The non-linear system of equation 3.30 is then discretised in time using a backward Euler's method and the non-linear system of equations, which is derived after the time discretisization, is solved using a Newton-Raphson procedure. Once the system is solved it is possible to retrieve the velocity and pressure fields, at each time step, using the values of the $\boldsymbol{a}$ anf $\boldsymbol{b}$ vectors and the fields approximation as presented in equation 3.1 and 3.2.

## 4. A NUMERICAL EXAMPLE

In the present section we will discuss the results obtained in the first application of the proposed model reduction procedure for Navier–Stokes flows. Given all the aforementioned features, the fluid dynamic problem considered is that of the vortex shedding caused by the low Re flow past a circular cylinder with main axis perpendicular to the undisturbed stream velocity $\mathbf{U}_\infty$. This is a well known and studied benchmark widely discussed and treated in literature from both the experimental and numerical point of view [27]. Due to the considerably larger extension of the cylinder in its axial direction and to the synchronisation observed in the vortices detaching from the cylinder at different axial locations, the vortex shedding mechanism exhibits an intrinsic two-dimensional nature. For such reason, the laminar Navier–Stokes equations (Equation (2.1)) written in a 2D domain represent a suitable model for the flow at hand. The resulting 2D computational grid is depicted in Figure 2, which also shows the domain height and width as a function of the cylinder diameter $D = 0.027$ m. The structured grid accounts for 13296 quadrilateral cells. The fluid considered in the simulations is water, having constant density $\rho_w = 1000 \text{kg/m}^3$ and kinematic viscosity $\nu = 10^{-6} \text{m}^2/\text{s}$. The boundary conditions are set according to Table 4. In each flow simulation, the fluid is started from rest and impulsively accelerated through the imposition of uniform and constant horizontal velocity $\mathbf{U}_\infty$ at the inlet boundary. Each simulation evolves in time until a final periodic regime solution is reached, and is then finally carried on for about 20 to 25 periods. The fluid dynamic drag and lift forces coefficients time history over the latter part of the simulation is then used to carry out Fourier analysis and assess the main vortex shedding frequency.

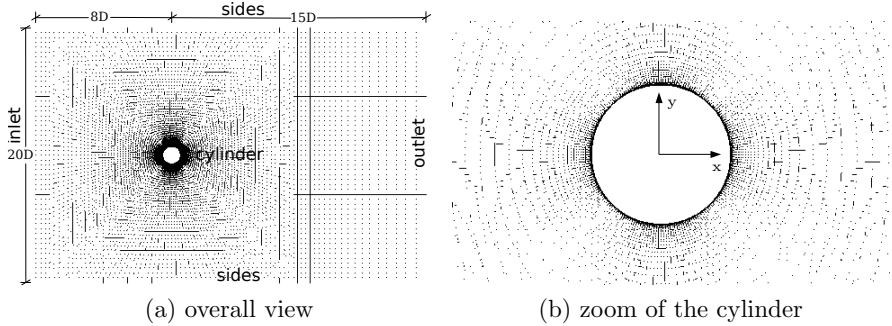

(a) overall view     (b) zoom of the cylinder

FIGURE 2. A sketch of the structured computational grid used for the high fidelity simulations. The picture also shows the main dimensions of the computational domain $\Omega_f$, as a function of the cylinder diameter $D = 0.027$ m.

|   | inlet | outlet | cylinder | sides |
|---|---|---|---|---|
| $\boldsymbol{u}$ | $\boldsymbol{u_{in}} = [u_{x_{in}}, 0]$ | $\nabla \boldsymbol{u} \cdot \boldsymbol{n} = \boldsymbol{0}$ | $\boldsymbol{u} = \boldsymbol{0}$ | $\boldsymbol{u} \cdot \boldsymbol{n} = 0$ |
| $p$ | $\nabla p \cdot \boldsymbol{n} = 0$ | $p = 0$ | $\nabla p \cdot \boldsymbol{n} = 0$ | $\nabla p \cdot \boldsymbol{n} = 0$ |

TABLE 1. Boundary Conditions

### 4.1. Constant Inflow Velocity.

In this first test case considered, an inlet velocity of $\boldsymbol{u_{in}} = [3.7\mathrm{e}{-3}, 0]$m/s, corresponding to $Re = 100$, is prescribed. The results of the HF simulations at several time instants are then used as the snapshots needed to compute the basis functions and set up the ROM model. Not only the comparison between HF and ROM solution provides a first assessment of the overall ROM performance, but it also helps understanding what is the effect of each value of the ROM parameter (that is the inlet velocity in our case) on the reduced solution. In particular, we will investigate on the number of bases required for a suitably accurate ROM solution. It will be investigated the response of the ROM also for a longer time window, larger respect to the one used to train the ROM model (Figure 8).

4.1.1. *Details of the full order simulation.*
The computational grid is the one presented in § 4. The convective term is discretised in space making use of the Gauss's theorem (see Equation (2.5)). The face center values of the variables are obtained from the center cell ones, which are the numerical problem unknowns, with an interpolation scheme consisting into a combination of a linear and upwind scheme. The diffusive term is discretised (see Equation (2.6)) in a similar fashion. In this case though, a central differencing interpolation scheme with non-orthogonality correction is preferred. Also the pressure gradient is discretised making use of Gauss's theorem (see equation (2.4)). Here, the face center pressure values are obtained from the cell center ones by means of a linear interpolation scheme, in which a limiting operation on the gradient is performed so as to preserve the monotonicity condition and ensure that the extrapolated face value is bounded by the neighbouring cell values [18]. As for the time discretisation, a backward Euler scheme is used. The overall time extent of the simulation is equal to $T = 3645$s, which is sufficiently long to reach a perfectly periodic response of the lift and drag forces. The simulation is run in parallel on 4 Intel® Core™ i7-4710HQ 2.50GHz processors, taking $T_{CPU_{HF}} = 1483$s $\approx$ 25min to be completed.

4.1.2. *Details of the ROM simulation.*
The ROM is constructed using the methodologies described in § 3. For the generation of the POD spaces, we considered 120 snapshots of the velocity, mass flux and pressure fields. The snapshots are collected in a time window covering approximately 1.5 periods of the vortex shedding phenomenon. More precisely, the last 73s of the HF simulation are used. The first two modes for velocity and pressure field respectively are presented in Figure 3 and 4. The ROM simulations are carried out using different values of the POD velocity space dimension $N_u = 3, 5, 7, 10$. The dimension of the POD pressure and mass flux space is set equal to the dimension of the velocity POD space $N_u = N_p$ but, for the way the ROM has been developed, also other choices are possible. The ROM simulation is run in serial, on the same processor used for the HF simulation. In this case, the time advancing of the ROM problem is carried out using a backward Euler's method. Reproducing the full 3645s extent of the high fidelity simulation requires, using the ROM model with the highest dimension of the POD space, approximately $T_{CPU_{ROM}} = 9.10s$. This corresponds to a speedup $SU \approx 650$.

4.1.3. *Analysis of the results.*
Using the settings described in the previous paragraph, four different ROM simulations are run, each featuring a different value of the POD space dimension. The results are compared with those of the HF simulation in terms of history of the lift and drag coefficients. In Figures 5 and 6 the time window used for the comparison is the same window used for the collection of the snapshots while in Figure 8 the comparison is performed on a different time window, larger respect to the one used to train ROM model. In Table 4.1.3 are reported the cumulative eigenvalues for the case with only one velocity and for the case with 5 different inlet velocities of subsection 4.2. The lift coefficient comparison is reported in Figure 5, while the drag coefficient time histories is presented in Figure 6. In each Figure, the four different plots refer to the four different values of the POD space dimensions considered in the ROM simulations. Finally, in Figure 7 the comparison is shown directly on the velocity and the pressure fields. In this case, the time step considered is the last one of the simulations, corresponding to $T = 3645$s. The left plot in Figure 7 refer to the velocity (top) and pressure (bottom) fields computed with the high fidelity simulations. The right plots refer to the velocity (top) and pressure (bottom) fields computed with the ROM, in which the POD space dimension has be set to $N_u = 10$. The plots show that, at a glance the HF and ROM solutions cannot be distinguished.

To provide a more quantitative evaluation of the error in the force coefficients reconstruction, for each ROM simulation we computed the Weighted Absolute Percentage Error (WAPE) [23] with respect to the HF simulations. For the case of the lift coefficient the WAPE has been directly applied to the lift signal without any modification, namely

$$(4.1) \qquad \varepsilon_{L_c} = \frac{100}{n} \sum_{t=1}^{n} \left| \frac{L_{c_t}^{HF} - L_{c_t}^{ROM}}{L_c^{HF}} \right| \%$$

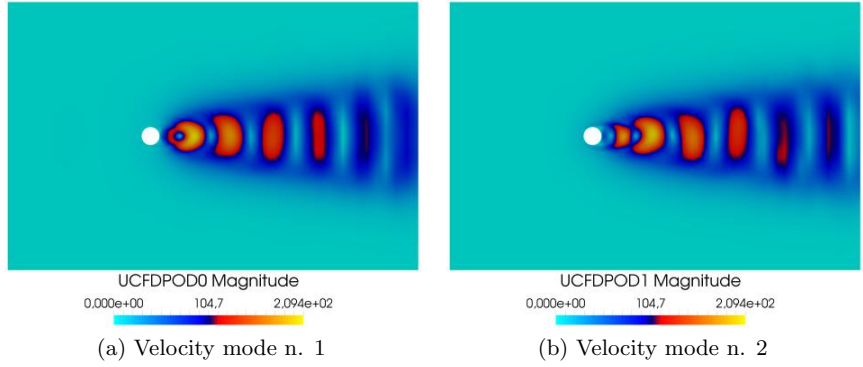

(a) Velocity mode n. 1  (b) Velocity mode n. 2

FIGURE 3. First two modes of the velocity field

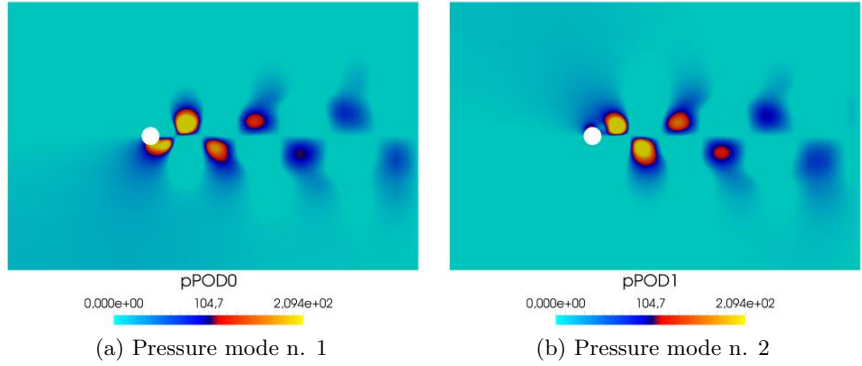

(a) Pressure mode n. 1  (b) Pressure mode n. 2

FIGURE 4. First two modes of the pressure field

|  | $N_u = N_p = 3$ | $N_u = N_p = 5$ | $N_u = N_p = 7$ | $N_u = N_p = 10$ |
|---|---|---|---|---|
| $\varepsilon_{L_c}(\%)$ | 11.50 | 4.32 | 1.59 | 1.89 |
| $\varepsilon_{D_c}(\%)$ | 64.49 | 13.94 | 4.69 | 6.43 |

TABLE 2. Error on $L_c$ and $D_c$ using different dimensions of the POD spaces

where $n_t$ is the number of sampling points, $L_{c_i}^{HF}$ and $L_{c_i}^{ROM}$ are the lift coefficients for the HF and ROM case respectively at the $i$−th time step. For the case of the drag coefficient, the WAPE has been applied to $D_c'$ that is the value of the drag shifted by the its mean value $D_c' = D_c - \overline{D_c}$. From Table 4.1.3 and Figures 5–6 one can see that, for the present case, adding more than 7 modes does not increase the accuracy of the ROM results. In Figure 8 the HF and ROM model have been simulated for a longer time window larger respect to the one used to train the ROM model. The time window used to calibrate the ROM model is indicated by the black arrow. Using only three modes the ROM model, for long time integration exhibits numerical instabilities as observed in [5]. Using more modes, including a percentage of energy up to 99%, instability phenomena vanish.

4.2. **Varying Inflow Velocity.**
In the second example the inlet velocity is used as a physical parameter and the HF simulations are performed using five different values of the *Reynolds* number $Re = [100, 125, 150, 175, 200]$. Using the same procedure described in § 4.1, 120 snapshots are collected for each different value of Re. Thus, the POD is performed on the resulting ensemble of 600 snapshots gathered from the 5 different HF simulations. According

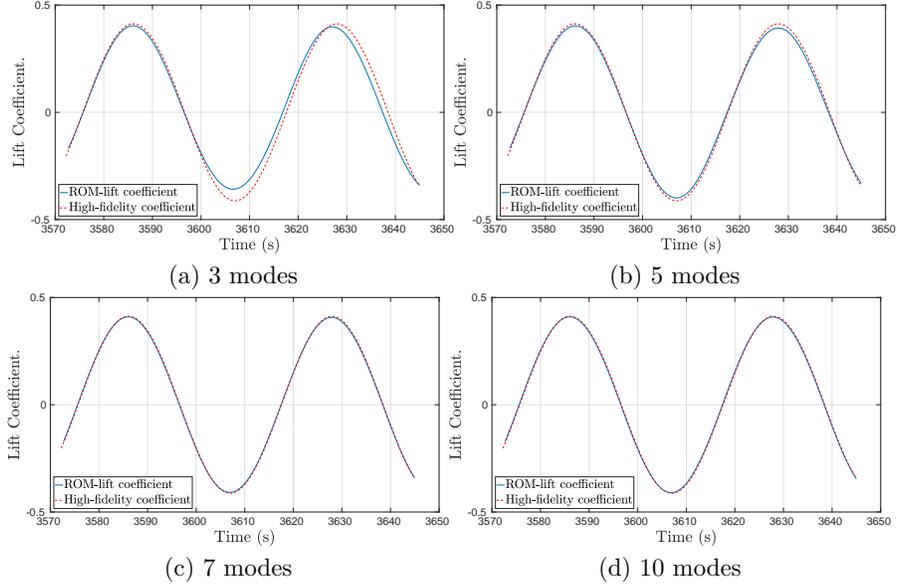

Figure 5. Comparison of the lift coefficient obtained with the HF and ROM simulations. The comparison is plotted on the same window used for the collection of the snapshots.

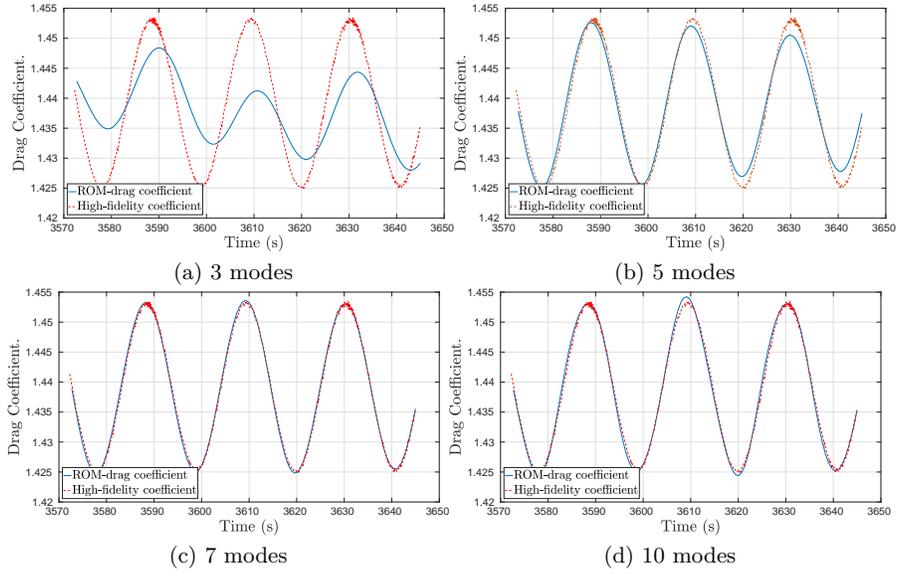

Figure 6. Comparison of the drag coefficient obtained with the HF and ROM simulations. The comparison is plotted on the same window used for the collection of the snapshots

to the results presented in § 4.1, 7 modes for each different velocity are considered In this case, leading to a total number of modes equal to 35. The computational details of the HF simulations are the same described in § 4.1. Once the offline phase is carried out and the ROM is set up, several reduced simulations are performed. Such simulations are also featuring inlet velocities that correspond to *Reynolds* numbers which were not considered in the HF analysis. Thus, this numerical experiment is devised to test if the ROM developed is able to reproduce the dependence of the system output with respect to an input parameter such as the velocity of the stream in which the cylinder is embedded. The comparison between the HF and ROM simulations is performed comparing the frequency corresponding to the peak of the power spectral density of the lift coefficient. Investigating such dependence is particularly important, as variations in the main stream velocity might result in different frequency components in the hydrodynamic force on the cylinder. By an engineering standpoint, it is typically very important to assess whether such frequency components are close to the structural natural frequency and might lead to resonance. The comparison can be seen in Figure 9 where the blue line with circles and the red line with stars refer respectively to the results of ROM and HF simulations. As one observes from the Figure the ROM simulations match well the HF simulations that for this range of *Reynolds* numbers lay on the red straight line.

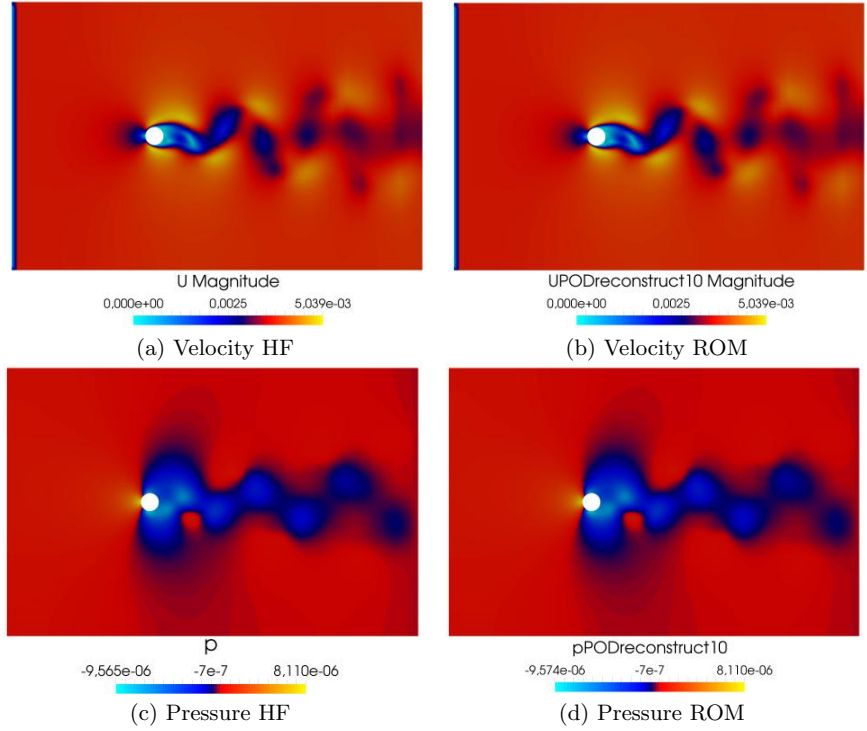

(a) Velocity HF  (b) Velocity ROM

(c) Pressure HF  (d) Pressure ROM

Figure 7. Comparison between velocity and pressure HF-ROM

| N Modes | $\boldsymbol{u}(\mathrm{Re} = 100)$ | $p(\mathrm{Re} = 100)$ | $\boldsymbol{u}(\mathrm{Re} = 100 : 200)$ | $p(\mathrm{Re} = 100 : 200)$ |
|---|---|---|---|---|
| 1 | 0.40509 | 0.906293 | 0.351995 | 0.840604 |
| 2 | 0.726335 | 0.951848 | 0.624987 | 0.916326 |
| 3 | 0.960515 | 0.992383 | 0.893356 | 0.971647 |
| 4 | 0.978525 | 0.996043 | 0.915563 | 0.981464 |
| 5 | 0.995695 | 0.999308 | 0.937393 | 0.987549 |
| 6 | 0.997671 | 0.999599 | 0.956756 | 0.993371 |
| 7 | 0.999518 | 0.999875 | 0.973212 | 0.995778 |
| 8 | 0.999732 | 0.99992 | 0.984389 | 0.997518 |
| 9 | 0.99994 | 0.999962 | 0.987669 | 0.998013 |
| 10 | 0.999967 | 0.999969 | 0.990883 | 0.998502 |

Table 3. Cumulative Eigenvalues

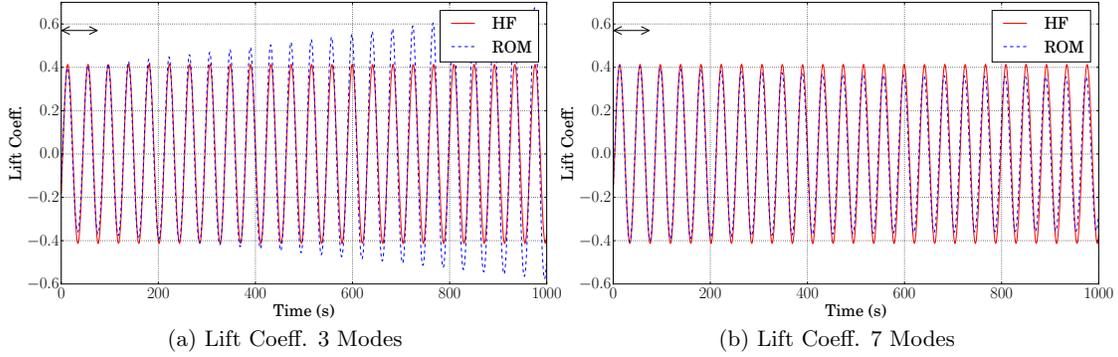

(a) Lift Coeff. 3 Modes

(b) Lift Coeff. 7 Modes

FIGURE 8. Comparison of lift coefficient obtained with HF model (continuous red line) and ROM model (dashed blue line). The left plot (a) refers to the ROM results obtained making use of 3 modes. The right plot (b) depicts the ROM results obtained with 10 modes. The black arrow on top left of the plots indicates the time window used to train the reduced order model.

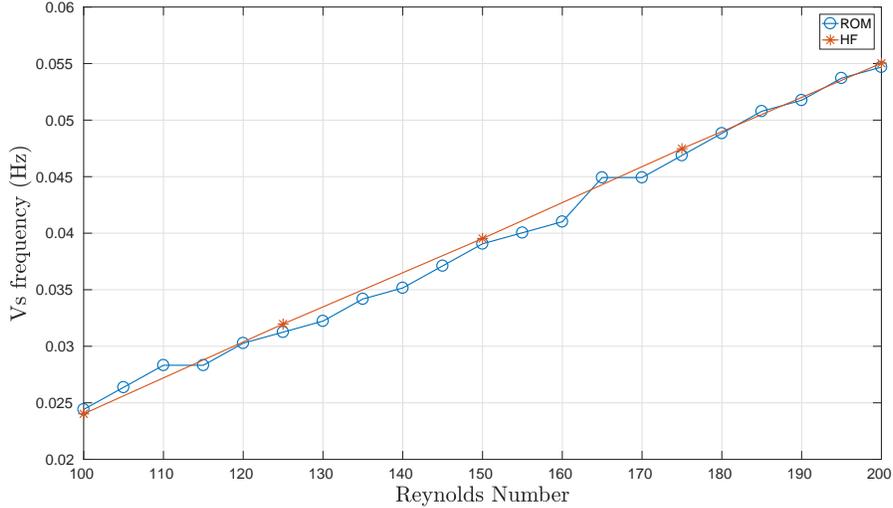

FIGURE 9. Comparison of ROM and HF for the case of increasing velocities - frequency of vortex shedding

## 5. Conclusions and future developments

In this work a POD-Galerkin ROM based on finite volume high fidelity simulations is presented. The ROM is generated such that to be fully consistent with the full order model, and both velocity and pressure fields are considered. In particular, the reconstruction of the reduced pressure field is carried out through the projection of the Poisson equation for pressure onto the POD pressure space. The ROM is then applied to approximate the unsteady viscous flow around a circular cylinder. In particular, the focus and originality is into reconstruction of the reduced pressure field using the Poisson pressure equation in a finite volume context. A reliable reconstruction of the pressure field permits the accurate reproduction of the lift and drag forces associated with the vortex shedding phenomenon. The ROM developed demonstrated to be capable of reproducing all the main features of the physical phenomenon in an accurate manner leading to a considerable computational time reduction (speedud $SU \approx 650$).

Also in the case of varying inflow velocities the ROM has demonstrated the ability of capturing the dependence of frequency of vortex shedding on the inflow velocity. As future developments the interest is into different efficient methodologies for the reconstruction of the pressure term and in particular to study the applicability of well known stabilization methods, used in the context of Galerkin ROM for finite elements [35, 3], to a finite volume framework [39]. The interest is also in the study of the same physical problem, but where the fluid-structure interaction problem is also considered. This will introduce additional complexities such as the mesh motion and additional equations to account for the structural dynamics and the fluid structure interaction coupling. Adding also the structural part to the ROM will be essential to tackle real-world engineering problems.

## Acknowledgements

We acknowledge the support provided by the European Research Council Consolidator Grant project Advanced Reduced Order Methods with Applications in Computational Fluid Dynamics - GA 681447, H2020-ERC COG 2015 AROMA-CFD and INdAM-GNCS.

## References

1. Imran Akhtar, Ali H. Nayfeh, and Calvin J. Ribbens, *On the stability and extension of reduced-order Galerkin models in incompressible flows*, Theoretical and Computational Fluid Dynamics **23** (2009), no. 3, 213–237.
2. J. Baiges, R. Codina, and S.R. Idelsohn, *Reduced-order modelling strategies for the finite element approximation of the incompressible Navier-Stokes equations*, Computational Methods in Applied Sciences **33** (2014), 189–216.
3. Francesco Ballarin, Andrea Manzoni, Alfio Quarteroni, and Gianluigi Rozza, *Supremizer stabilization of POD-Galerkin approximation of parametrized steady incompressible Navier–Stokes equations*, International Journal for Numerical Methods in Engineering **102** (2015), no. 5, 1136–1161.
4. Maxime Barrault, Yvon Maday, Ngoc Cuong Nguyen, and Anthony T. Patera, *An 'empirical interpolation' method: application to efficient reduced-basis discretization of partial differential equations*, Comptes Rendus Mathematique **339** (2004), no. 9, 667 – 672.
5. M. Bergmann, C.-H. Bruneau, and A. Iollo, *Enablers for robust POD models*, Journal of Computational Physics **228** (2009), no. 2, 516–538.
6. J. Burkardt, M. Gunzburger, and H.-C. Lee, *POD and CVT-based reduced-order modeling of Navier-Stokes flows*, Computer Methods in Applied Mechanics and Engineering **196** (2006), no. 1-3, 337–355.
7. Alfonso Caiazzo, Traian Iliescu, Volker John, and Swetlana Schyschlowa, *A numerical investigation of velocity-pressure reduced order models for incompressible flows*, Journal of Computational Physics **259** (2014), 598 – 616.
8. Kevin Carlberg, Charbel Farhat, Julien Cortial, and David Amsallem, *The GNAT method for nonlinear model reduction: Effective implementation and application to computational fluid dynamics and turbulent flows*, Journal of Computational Physics **242** (2013), 623 – 647.
9. F. Chinesta, A. Huerta, G. Rozza, and K. Willcox, *Model Order Reduction*, Encyclopedia of Computational Mechanics (In Press, 2017).
10. Francisco Chinesta, Pierre Ladeveze, and Elías Cueto, *A Short Review on Model Order Reduction Based on Proper Generalized Decomposition*, Archives of Computational Methods in Engineering **18** (2011), no. 4, 395.
11. A. E. Deane, I. G. Kevrekidis, G. E. Karniadakis, and S. A. Orszag, *Low-dimensional models for complex geometry flows: Application to grooved channels and circular cylinders*, Physics of Fluids A: Fluid Dynamics **3** (1991), no. 10, 2337–2354.
12. A. Dumon, C. Allery, and A. Ammar, *Proper general decomposition (PGD) for the resolution of Navier-Stokes equations*, Journal of Computational Physics **230** (2011), no. 4, 1387–1407.
13. M.L. Facchinetti, E. de Langre, and F. Biolley, *Coupling of structure and wake oscillators in vortex-induced vibrations*, Journal of Fluids and Structures **19** (2004), no. 2, 123 – 140.

14. W. R. Graham, J. Peraire, and K. Y. Tang, *Optimal control of vortex shedding using low-order models. Part I:open-loop model development*, International Journal for Numerical Methods in Engineering **44** (1999), no. 7, 945–972.
15. Ronald T Hartlen and Iain G Currie, *Lift-oscillator model of vortex-induced vibration*, Journal of the Engineering Mechanics Division **96** (1970), no. 5, 577–591.
16. Jan S Hesthaven, Gianluigi Rozza, and Benjamin Stamm, *Certified Reduced Basis Methods for Parametrized Partial Differential Equations*, Springer International Publishing, 2016.
17. R.I. Issa, *Solution of the implicitly discretised fluid flow equations by operator-splitting*, Journal of Computational Physics **62** (1986), no. 1, 40–65.
18. Hrvoje Jasak, *Error analysis and estimation for the finite volume method with applications to fluid flows*, Ph.D. thesis, Imperial College, University of London, 1996.
19. I. Kalashnikova and M. F. Barone, *On the stability and convergence of a Galerkin reduced order model (ROM) of compressible flow with solid wall and far-field boundary treatment*, International Journal for Numerical Methods in Engineering **83** (2010), no. 10, 1345–1375.
20. K. Kunisch and S. Volkwein, *Galerkin proper orthogonal decomposition methods for a general equation in fluid dynamics*, SIAM Journal on Numerical Analysis **40** (2002), no. 2, 492–515.
21. Stefano Lorenzi, Antonio Cammi, Lelio Luzzi, and Gianluigi Rozza, *POD-Galerkin method for finite volume approximation of Navier-Stokes and RANS equations*, Computer Methods in Applied Mechanics and Engineering **311** (2016), 151 – 179.
22. X. Ma and G.E. Karniadakis, *A low-dimensional model for simulating three-dimensional cylinder flow*, Journal of Fluid Mechanics **458** (2002), 181–190.
23. Spyros Makridakis, *Accuracy measures: theoretical and practical concerns*, International Journal of Forecasting **9** (1993), no. 4, 527 – 529.
24. F. Moukalled, L. Mangani, and M. Darwish, *The Finite Volume Method in Computational Fluid Dynamics: An Advanced Introduction with OpenFOAM and Matlab*, 1st ed., Springer Publishing Company, Incorporated, 2015.
25. Bernd R. Noack and Helmut Eckelmann, *A low-dimensional Galerkin method for the three-dimensional flow around a circular cylinder*, Physics of Fluids **6** (1994), no. 1, 124–143.
26. Bernd R. Noack, Paul Papas, and Peter A. Monkewitz, *The need for a pressure-term representation in empirical Galerkin models of incompressible shear flows*, Journal of Fluid Mechanics **523** (2005), 339–365.
27. Michael P. Païdoussis, *Fluid-Structure Interactions. Slender Structures and Axial Flow. Volume 1.*, first ed., Academic Press, 1998.
28. ______, *Fluid-Structure Interactions. Slender Structures and Axial Flow. Volume 2.*, first ed., Academic Press, 2003.
29. S.V Patankar and D.B Spalding, *A calculation procedure for heat, mass and momentum transfer in three-dimensional parabolic flows*, International Journal of Heat and Mass Transfer **15** (1972), no. 10, 1787 – 1806.
30. Alfio Quarteroni, Andrea Manzoni, and Federico Negri, *Reduced Basis Methods for Partial Differential Equations*, Springer International Publishing, 2016.
31. Alfio Quarteroni and Gianluigi Rozza, *Numerical solution of parametrized Navier–Stokes equations by reduced basis methods*, Numerical Methods for Partial Differential Equations **23** (2007), no. 4, 923–948.
32. G. Rozza, D. B. P. Huynh, and A. Manzoni, *Reduced basis approximation and a posteriori error estimation for Stokes flows in parametrized geometries: Roles of the inf-sup stability constants*, Numerische Mathematik **125** (2013), no. 1, 115–152 (English).
33. G.a Rozza, D.B.P.b Huynh, and A.T.c Patera, *Reduced basis approximation and a posteriori error estimation for affinely parametrized elliptic coercive partial differential equations: Application to transport and continuum mechanics*, Archives of Computational Methods in Engineering **15** (2008), no. 3, 229–275.
34. Gianluigi Rozza, *Reduced basis methods for Stokes equations in domains with non-affine parameter dependence*, Computing and Visualization in Science **12** (2009), no. 1, 23–35.
35. Gianluigi Rozza and Karen Veroy, *On the stability of the reduced basis method for Stokes equations in parametrized domains*, Computer Methods in Applied Mechanics and Engineering **196** (2007), no. 7, 1244 – 1260.
36. S. Sirisup and G.E. Karniadakis, *Stability and accuracy of periodic flow solutions obtained by a POD-penalty method*, Physica D: Nonlinear Phenomena **202** (2005), no. 3-4, 218 – 237.

37. Lawrence Sirovich, *Turbulence and the Dynamics of Coherent Structures part I: Coherent Structures*, Quarterly of Applied Mathematics **45** (1987), no. 3, 561–571.
38. Giovanni Stabile, Hermann G Matthies, and Claudio Borri, *A novel reduced order model for vortex induced vibrations of long flexible cylinders*, Submitted to Journal of Ocean Engineering (2016).
39. Giovanni Stabile and Gianluigi Rozza, *Stabilized Reduced order POD-Galerkin techniques for finite volume approximation of the parametrized Navier–Stokes equations*, submitted (2017).
40. Vincenz Strouhal, *Über eine besondere Art der Tonerregung*, Annalen der Physik **241** (1878), no. 10, 216–251.
41. H. K. Versteeg and W. Malalasekera, *An Introduction to Computational Fluid Dynamics. The Finite Volume Method*, Longman Group Ltd., London, 1995.
42. Henry G Weller, G Tabor, Hrvoje Jasak, and C Fureby, *A tensorial approach to computational continuum mechanics using object-oriented techniques*, Computers in physics **12** (1998), no. 6, 620–631.
43. D. Xiao, F. Fang, A.G. Buchan, C.C. Pain, I.M. Navon, J. Du, and G. Hu, *Non linear model reduction for the navier stokes equations using residual deim method*, Journal of Computational Physics **263** (2014), 1 – 18.
44. Momchilo M Zdravkovich, *Flow around circular cylinders: Volume 1: Fundamentals*, vol. 350, Cambridge University Press, 1997.
45. ———, *Flow around circular cylinders: Volume 2: Applications*, vol. 2, Oxford University Press, 2003.